\documentclass[10pt]{article}
\usepackage[english]{babel}
\usepackage{latexsym}
\usepackage{color}
\usepackage{amsmath}
\usepackage[latin1]{inputenc}
\usepackage{amsfonts}
\usepackage{amssymb}
\usepackage{mathrsfs}
\usepackage{eucal}
\usepackage{exscale}
\usepackage{makeidx}
\usepackage{makeidx}
\usepackage{graphicx}
\usepackage[T1]{fontenc}
\usepackage{amsmath,amssymb}
\usepackage[all]{xy}
\newcommand{\qed}{{$\square$}\medbreak}

\newtheorem{thm}{Theorem}
\newtheorem{rem}{Remark}
\newtheorem{df}{Definition}

\newtheorem{cor}{Corollary}
\newtheorem{lem}{Lemma}

\newtheorem{exemples}{Examples}
\newtheorem{exemple}{Example}

\newcommand{\preuve}{\indent {\it Proof.}\hspace{4mm}}

\title{h-Amalgamation Bases in The Class of Non Trivial Abelian 
Groups }
\author{Mohammed Belkasmi}

\begin{document}
\maketitle

\begin{abstract} 
In this paper we give a complete description of the h-amalgamation
bases in the class of non trivial abelien groups.
 \end{abstract} 
\section*{Introduction}
The amalgamation property is essentially a property of certain types of structures and their homomorphisms. It ensures the existence of common extensions for structures sharing certain substructures. The concept has been adapted and applied in various mathematical contexts to study the structure and behavior of different algebraic and model-theoretic systems.

The amalgamation property in group theory is associated with the study of free products and amalgamated products of groups. 
The concept of an amalgamated free product of groups and the associated amalgamation property has evolved over time. In 1926, Philip Hall introduced this concept, which involves amalgamating two subgroups over a common subgroup. The development of these concepts has been the subject of various mathematical studies and proofs, such as Kurosh's proof in 1934 that a subgroup of a free product of groups is again a free product \cite{5}.
The amalgamated free product has been the subject of various mathematical investigations and attempts to generalize its 
properties \cite{2}. Additionally, research has been conducted on the free product of two groups with an amalgamated subgroup of finite index in each factor \cite{3}. These studies have contributed to the ongoing development and understanding of the amalgamation property in group theory.
The amalgamation property continues to be an active area of research within group theory. Mathematicians explore its connections with other topics in algebra, geometry, and topology.\\

The algebraic meaning of the amalgamation property lies in its ability to combine smaller algebraic structures into larger ones while maintaining specific relationships or properties. This property is essential for studying and understanding the interactions and transformations within algebraic systems, providing a powerful tool for algebraic analysis and exploration (\cite{lattice}).\\
In model theory, the amalgamation property plays a crucial role in the study of theories and their models. The presence or absence of the amalgamation property can help classify theories based on their behavior. It is closely related to concepts such as 
the quantifiers elimination property (\cite{marker}), 
stability, independence, and categoricity of theories (\cite{stable}).

The study of the amalgamation property in algebra and model theory heavily relies on embeddings. However, the positive logic offers a more comprehensive framework that encompasses various forms of amalgamation, such as asymmetric amalgamation and h-amalgamation. For further details,  refer to (\cite{ana2}). 
The particularity of positive model theory lies essentially in the study of h-inductive theories, and it prohibits the use of the negation operator in formulas construction. Thus, it focuses on positive formulas rather than general formulas, and on homomorphisms rather than embeddings. This approach creates new situations in the study of existentially closed structures and the amalgamation property, which go beyond the traditional framework of first-order logic.  For example, but not exclusively, in the framework of positive logic, the class of existentially closed groups is reduced to the trivial group, and each group has the property of amalgamation. However, in the framework of  the first order the class logic,  of existentially closed groups is not axiomatizable. So, to prevent such undesirable implications, modifications are made to the language of the theory in order to eliminate the trivial structures from the class of models associated with the theory.

This paper focuses specifically on the exploration of the class of  h-amalgamation bases  of the theory of non trivial abelian groups 
in the context of positive logic. The organization of the article is as follows: In Section 1, we establish the foundation of the essential tools of positive logic.  Section 2 is devoted to the investigation of the h-amalgamation bases of the class of non trivial abelien groups. We show that  the h-amalgamation bases are the groups 
that contain only one subgroup isomorphic to $\mathbb{Z}/p^n\mathbb{Z}$.
\section{Positive theory of non trivial abelian groups}
In this section, we provide an overview of the fundamental concepts of positive mathematical logic within the context of the theory of abelian groups, emphasizing the crucial concepts  for this paper.

Consider a first-order language $L$, consisting of a set of symbols of functions, symbols of relations and symbols of constants. An $L$-structure is defined as a set $A$ that meets the following conditions:
\begin{itemize}
\item Every symbol of n-ary function $f$ of 
$L$ is interpreted by a function $f$ defined from $A^n$ to $A$.
\item Every m-ary relation symbol $R$ in the language $L$ has a corresponding interpretation in the structure $A$, represented by a subset of $A^m$.
\item Every symbol of constant of $L$ has an  interpretation  by an
element of $A$.
\end{itemize}
Let $A$ and $B$ be two $L$-structures. A function $f$ from $A$ to $B$ is called a $L$-homomorphism if it  satisfies  the following properties:
\begin{itemize}
\item For every symbol of n-ary function $h$ of $L$, and 
for every $\bar a=(a_1,\cdots,a_n)\in A^n$; $f(h(\bar a))=h(\overline{f(a_i)})$, where $\overline{f(a_i)}=(f(a_1),
\cdots, f(a_n))$.
\item For every symbol of n-ary relation $R$ of $L$, and 
for every $\bar a\in A^n$; if $A$ satisfies $R(\bar a)$ then
$B$ satisfies $R(\overline{f(a_i)})$.
\item For every symbols of constant $c$ in $L$, we have 
$f(c_A)=c_B$.
\end{itemize}
The h-inductive $L$-sentences are  build by the conjunctions of expressions
 of the form:
$$\forall\bar y, \exists\bar a;\,\, \varphi(\bar x, \bar y)
\rightarrow\psi(\bar x, \bar y)$$
where $\varphi(\bar x, \bar y)$ and $\psi(\bar x, \bar y)$
are existential $L$-formulas (see Example \ref{exm illustr}).\\
Ah-inductive $L$-theory $T$ is  a set of h-inductive $L$-sentences 
 that can be satisfied by a $L$-structure. Every $L$-structure that satisfies
a h-inductive theory $T$ is said a model of $T$. In the rest of the paper,  we will employ the terms "L-theory" and "L-sentences" instead of "h-inductive L-theory" and "h-inductive L-sentences," respectively.\\
For further details, \cite{ana} is a sufficiently complete reference.

\begin{exemples}
\item \textbf{Groups}:
 The language of groups is given by the set  $L_g = \{e, \cdot, ^{-1}\}$, 
where $e$ is a symbol of constant, $\cdot$, a symbol of function
of arity 2, and, $^{-1}$, a symbol of function of arity $1$.\\
A positive $L_g$-formula $\varphi(\bar x, \bar y)$ is
finite conjunctions and disjunctions  of formulas 
of the form:
$$\exists\bar y;\ \ \ x_{\sigma(1)}^{n_1}\cdot y_{\delta(1)}^{n_1}
\cdot x_{\sigma(2)}^{n_2}\cdot y_{\delta(2)}^{n_2}\cdots
x_{\sigma(p)}^{n_p}\cdot y_{\delta(q)}^{n_q}=e$$
where $\bar x=(x_1, x_2,\cdots, x_n)$, 
$\bar y=(y_1, y_2,\cdots, y_m)$ are variables, $\sigma$ 
(resp. $\delta$) runs 
overs the set of mapping defined from a finite subset of 
$\mathbb{N}$ into the set $\{1,\cdots, n\}$ 
(resp, $\{1,\cdots, m\}$).\\
 A group $G$ is a $L_g$-structure that satisfies the $L_g$-theory $T_g$,
the following  set of  $L_g$-sentences:
 \begin{itemize}
 \item $\forall x, y, z;\ x\cdot(y\cdot z)=(x\cdot y)\cdot z$.
 \item $\forall x;\ \  e\cdot x= x\cdot e= x.$
 \item $\forall x;\ \  x\cdot x^{-1}=  x^{-1}\cdot x=e.$
 \end{itemize}
where $x, y, z$ are symbols of variables. Note that, in the 
expression of  every $L$-sentence, the number of variables is 
equal to the number of variables that follow the universal quantifier $\forall$ and the existential quantifier $\exists$.\\
The $L_g$-homomorphisms are the homomorphisms of groups in 
the usual algebraic sense.
\item \textbf{Abelian groups}:
The language of abelian group is $L_{ab}=\{0, +, -\}$,  and
abelian groups are the models of $T_{ab}$, were $T_{ab}$
is the set of the following 
$L_{ab}$-sentences:
\begin{itemize}
\item $\forall x, y, z;\ x+(y+ z)=(x+ y)+ z$,
\item $\forall x, y;\ x+y=y+x$, 
\item $\forall x;\ x+0=x$,
\item $\forall x;\ x+(-x)=0$.
\end{itemize}
\item \textbf{Non trivial abelian groups}:
Let $L^*= L_{ab}\cup \{g\}$ be
the language of non trivial  abelian group, 
where $g$  is a symbol of constant, and $L_{ab}$ the 
language of abelian groups. Let 
$$T^*_{ab}=T_{ab}\cup \{g\neq 0\}$$ 
For every abelien group $G$ and for every $g\in G$, the pair 
$(G, g)$ is a model of $T^*_{ab}$ if and only if $g\neq 0$.\\
Let $(G, g)$ and $(K, k)$ be two models of $T^*_{ab}$ and 
$f$ a homomorphism  of groups defined from  $G$ to $K$.  $f$
is a $L^*$-homomorphism if and 
only if  $f(g)=k$.
\end{exemples}

 \section{h-Amalgamation}

 \begin{df}
A model $A$ of a $L$-theory $T$ is said to be an 
$h$-amalgamation basis of $T$ if for every   models $B$ and $C$ of $T$, $f$ a $L$-homomorphism from $A$ to $B$, and $g$ a $L$-homomorphism from $A$ to $C$, there 
are a model $D$ of $T$ and $L$-homomorphisms $f'$ and $g'$ such 
that the following diagram commutes:
\[
\xymatrix{
    A \ar[r]^{g} \ar[d]_{f} & {C} \ar[d]^{g'} \\
    B \ar[r]_{f'} & {D}
  }
\]
We say that $T$ has the $h$-amalgamation property if every model 
of $T$ is an $h$-amalgamation basis of $T$.
 \end{df}
\begin{rem}
Let $L$ be a language with at most one symbol of constant. Let 
$T$ be a $L$-theory and $A$ a model of $T$ that satisfy the 
following properties:
\begin{itemize}
\item $A=\{c\}$ is a model of $T$, where $c$ is the interpretation of the constant
of the language if the language contains a symbol of constant, otherwise, $A$ is a singleton model of $T$ by hypothesis.
\item For every model $B$ of $T$,  the constant 
mapping defined from $B$ into $A$  is a $L$-homomorphism.
\end{itemize}
$T$ has the h-amalgamation property. 
Indeed, for every $B, C$ and $D$ models of $T$, and for every 
$L$-homomorphisms $f$ and $g$ defined from $B$ respectively to 
$C$ and $D$, the following diagram commutes:
\[
\xymatrix{
    B \ar[r]^{f} \ar[d]_{g} & {C} \ar[d]^{g'} \\
    D \ar[r]_{f'} & {A}
  }
\]
where $f'$ and $g'$ are the constant homomorphisms.\\
Note that this type of amalgamation is not desirable, indeed,  the absorbing model $A$  destroy  the properties of the other models
of the theory. This phenomenon is observed in many theories, notably the theory of ordered set,  the theory of lattices, the theory of groups, etc. To address this issue, it is necessary to eliminate the absorbent model of theory, this can be achieved by modifying both the language and the theory, and this is precisely what we do in the definition of the language and the theory of non trivial abelian  groups,(see Example \ref{exm illustr}). 
\end{rem}

\begin{lem}\label{char ab gr amal bs}
Let $(G,g),  (K, k)$ and $(L, l)$ be three models of $T_{ab}^*$.
Let $f$ be a $L^*$-homomorphism from  $(G,g)$ to  $(K, k)$ and 
$h$ a $L^*$-homomorphism from  $(G,g)$ to  $(L, l)$. The following 
assertions are equivalent:
\begin{enumerate}
\item there are $(D, d)$ a model of $T_{ab}^*$, $f'$  a $L^*$-homomorphism from $(K, k)$  to  $(D, d)$, and 
$h'$ a $L^*$-homomorphism from  $(L, l)$ to  $(D, d)$, such
that $  f'\circ f(x)= h'\circ h(x)$ for every $x\in G$. 
\item $l\not\in h( ker(f))$ and $k\not\in f(ker( h))$.
\end{enumerate}
\end{lem}
\preuve
\begin{itemize}
\item $1\Rightarrow 2$:\\
Suppose that $k=f(a)$ for some $a\in ker( h)$. By 
the hypothesis of the assertion (1), we have 
 $f'\circ f(a)= h'\circ h(a)$. so 
 $$d=f'\circ f(a)= h'\circ h(a)=0,$$ 
then $d=0$, contradiction.
\item $2\Rightarrow 1$:\\
Suppose that $l\not\in h( ker (f))$ and $k\not\in f(ker (h))$. 
Consider the following subset of $K\times L$:
$$H=\{(f(a), h(-a))|\ a\in G\}$$
It is clear that $H$ is a subgroup of $K\times L$. 
 Let the following diagram:
\[
\xymatrix{
    (G, g) \ar[r]^{f} \ar[d]_{h} & {(K,k)} \ar[d]^{f'} \\
    (L,l) \ar[r]_{h'} & {K\times L/H}
  }
\]
where  $K\times L/H$ is the abelian quotient group $K\times L$ of $H$,
 $f'$ and $h'$ are the homomorphisms defined by 
$f'(x)=\overline{(x,0)}$  and 
$h'(y)=\overline{(0, y)}$. It is clear that the diagram commutes.
To finish the proof, it suffices to show that $f'$ and $h'$ are 
$L^*$-homomorphisms. Suppose that $f'$ is not a $L^*$-homomorphisms, then $f'(k)=\overline{(k,0)}= \overline{(0, 0)}$.
So $(k,0)\in H$, which implies that $k=f(a)$ and $h(a)=0$.
Thus $k\in f(ker(h))$, a contradiction.\qed
\end{itemize}
\begin{lem}\label{char amal}
An abelian group $(G, g)$ is an $h$-amalgamation basis of $T_{ab}^*$
if and only if it satisfies the following property:\\
For every proper subgroups $H$ and $K$ of $G$;
$g\in H+K$ if and only if  $g\in H$ or  $g\in K$
\end{lem}
\preuve 
Let $(G, g)$ be an $h$-amalgamation basis of $T_{ab}^*$. Let 
$H$ and $K$ be non trivial subgroups of $G$
such that  $g\not\in H\cup K$. Then the natural mappings 
$\pi_1, \pi_2$  
from $G$ into $G/H$ and $G/K$ respectively are $L^*$-homomorphisms. By Lemma 
\ref{char ab gr amal bs}, we have:
$$\pi_1(g)\not\in \pi_1(K),\ \  \pi_2(g)\not\in \pi_2(H)$$
which implies the following 
$$
\left\{
    \begin{array}{ll}
        \forall k\in K, & g-k\not\in H \\
         \forall h\in H, & g-h\not\in K
    \end{array}
\right.
$$
Thereby $g\not\in H+K$.\\
Conversely, assume that for every non trivial subgroups $H$ and 
$K$ of $G$, if $g\not\in H\cup K$ then $g\not\in H+ K$. Let 
the following schemas:
$$
\xymatrix{
(F,l)&\ar[l]_{h} (G, g)\ar[r]^{f}& (E, e)\\
}
$$
where $h$ and $f$  are $L^*$-homomorphisms.
Let $H=ker(f)$ and $K= ker(h)$, then $g\not\in H\cup K$, 
 so $g\not\in H+ K$. Suppose that $e\in f(K)$, let $e=f(k)$ 
 for some $k\in K$. So, $f(g)=f(k)$, which implies $g-k\in H$, 
 then $g\in H+K$, contradiction. By  Lemma \ref{char ab gr amal bs}, the schemas is h-amalgamable.\qed
\begin{cor}\label{corchar amal} 
Let $(G, g)$ be an $h$-amalgamation basis 
of $T_{ab}^*$ then $o(g)$ is a power of a prime number.
\end{cor}
\preuve 
Assume that $(G, g)$ is an $h$-amalgamation basis of $T_{ab}^*$. 
We claim that  $o(g)$ is finite. Indeed,
suppose that $o(g)$ is infinite, let $p$ and $q$ two co-prime integers, so there are $u, v\in\mathbb{Z}$ such that $1=up+vq$.
Then 
$$g=upg+vqg.$$
By Lemma \ref{char amal}, $g\in \langle pg\rangle$ or 
$g\in \langle qg\rangle$, then  $o(g)$
is finite.\\
Now, suppose that  $o(g)=pq$ where $p$ and $q$ are co-prime numbers. One repeat the same argument as above we obtain 
$o(g)$ divides $p$ or $q$, contradiction. Thereby $o(g)$
is the power of  a prime number.\qed
\begin{exemple}\label{exm illustr}
\begin{itemize}
\item For every prime number $p$ and for every integer $n$,
the group $(\mathbb{Z}/p^n\mathbb{Z}, a)$ where  $a\neq 0$
is an $h$-amalgamation basis of $T_{ab}^*$. Indeed, this results from the fact that  
the set of proper subgroups of $\mathbb{Z}/p^n\mathbb{Z}$ is
 totally ordered by inclusion.
 \item For every prime number $p$, the group
 $((\mathbb{Z}/p\mathbb{Z})^2, (a, b))$ 
 in not an $h$-amalgamation basis of $T_{ab}^*$. Indeed,
for every $(a, b)\in (\mathbb{Z}/p\mathbb{Z})^2$ we 
have 
\begin{equation} \label{sums}
 \left\{
    \begin{array}{lll}
        (a, b)\in <(1, 0)>+<(0,1)> & \mbox{if}\  a\neq 0\ 
       \mbox{and}\  b\neq 0\\
       (a, b)\in <(1, 1)>+<(0,1)>  & \mbox{if}\ b=0\\
       (a, b)\in <(1, 1)>+<(1,0)>  & \mbox{if}\ a=0.\\
    \end{array}
\right.
\end{equation}
However  $(a, b)$ does not belong to any  subgroups of the sums in the formulas (\ref{sums}). More generally, for every prime number
 $p$, for every integer $n\geq 2$, and for every 
 $(a_1,\cdots, a_n)\in (\mathbb{Z}/p\mathbb{Z})^n$, the group 
 $((\mathbb{Z}/p\mathbb{Z})^n, (a_1,\cdots, a_n))$ is not 
 an $h$-amalgamation basis of  $T_{ab}^*$.
\end{itemize}
\end{exemple}
\begin{thm}\label{thm char aml by uniquess}
Let $G$ be a non trivial abelian group and $p$ a prime number. 
If $G$ has only one maximal $p$-subgroup 
$H\simeq \mathbb{Z}/p^k\mathbb{Z}$ where $p$ is a prime and $k$
integer,  then for every $h\in H-\{0\}$, 
$(G, h)$ is an $h$-amalgamation basis of $T_{ab}^*$.
\end{thm}
\preuve 
Assume that  $H\simeq \mathbb{Z}/p^k\mathbb{Z}$ is  the unique
 p-subgroup of $G$, and let $h\in H-\{0\}$. Assume that 
 $h\in L+ K$, where $L$ and $K$ are two proper subgroups of $G$. 
 Without loss of generality we can replace the subgroups $L$ and 
 $K$ by the following subgroups:
 $$
\left\{
    \begin{array}{ll}
        L_g=\{l\in L|\ \exists k\in K,\ l+k\in<g>\} \\
        K_g=\{k\in K|\ \exists l\in L,\ l+k\in<g>\} 
    \end{array}
\right.
$$
We distinguish two cases.\\
Case 1. Assume  $ L_g\cap  K_g=\{0\}$. 
In this case  $ L_g$ and $ K_g$ are finite and $ |L_g+ K_g|=| L_g|\cdot|K_g|$.
Given that $<g>\leq L_g+ K_g $, then $p$ divides $| L_g|$ or 
$|K_g|$. Assume that both $|L_g|$ and $|K_g|$ contain  subgroups of order $p$.
By the uniqueness of $H$ in $G$ and the fact that 
$H\simeq \mathbb{Z}/p^k\mathbb{Z}$ contains a unique subgroup of 
order $p$, we have $L_g\cap K_g\neq\{0\}$, a contradiction.
Thereby $<g>\leq H_g$ or $<g>\leq K_g $.\\
Case 2. Assume that  $ L_g\cap  K_g\neq\{0\}$. Let $N= L_g\cap  K_g$, we have 
\begin{equation}\label{eq 1}
(L_g+K_g)/ N=L_g/N\oplus K_g/N.
\end{equation}
Suppose that $g\not\in N$ and $o(g)=p^n$. 
By definition of $L_g$,
for every  $h\in L_g$ there are $k\in K_g$ and $t\in\mathbb{N}$
 such that 
$h+k=t\cdot g$,
so $p^{n}h+p^{n}k=tp^n g=0$. Thus 
$p^{n}h, p^{n}k\in N$, which implies that 
$o(\bar h)$ in $L_g/N$ divides $p^{n}$. Thereby 
$L_g/N$ and $K_g/N$ are $p$-groups.\\
Now, we will show that $L_g/N$ and $K_g/N$ are cyclic groups.
 Let $\bar g=\bar h+\bar k$
where $h\in L_g$ and $k\in K_g$. For every $h'\in L_g$ there 
exist $k'\in K_g$ and $t\in\mathbb{N}$ such that
 $h'+k'=t g$. Given that 
 $t\bar g=t\bar h+t\bar k$ and 
 $\bar h'+\bar k'=t \bar g$, by the property
  (\ref{eq 1}) we obtain:
$$
 \left\{
    \begin{array}{ll}
       \bar h'= t\bar h \\
       \bar k'= t\bar k
    \end{array}
\right.
$$ 
Which implies that  $L_g/N$ and $K_g/N$ are cyclic and generated
respectively by $\bar h$ and $\bar k$.\\
On the other hand,  suppose that $o(\bar g)=p^m$ in $G/N$. Since
$\bar g=\bar h+\bar k$, then 
$o(\bar g)=lcm(o(\bar h), o(\bar k))$ (because
$L_g/N\cap K_g/N=\{\bar 0\}$). Given that $L_g/N$ and $K_g/N$
are cyclic p-groups, then $o(\bar h), o(\bar k)$ are powers of $p$. 
So $o(\bar g)=max(o(\bar h), o(\bar k))$. Thereby 
$o(\bar g)=o(\bar h)$ or $ o(\bar g)=o(\bar k)$. Suppose that 
$o(\bar g)=o(\bar h)$. Then
\begin{equation}\label{eq2}
\vert L_g/N\vert=o(\bar g)=p^m.
\end{equation} 
Since $p^m g\in N$ and 
$o(p^m g)=p^{n-m}$ in $N$ (ie. $p^n$ is the order of $g$ in $G$), then $p^{n-m}$ divides
$|N|$. Thereby, from (\ref{eq2}), $p^n$ divides $|L_g|$.
By the uniqueness of subgroup of order 
$p^n$ in $H\simeq \mathbb{Z}/p^k\mathbb{Z}$  we obtain  $g\in L_g$.\qed
\begin{lem}
A non trivial abelian group $(G, g)$ is an $h$-amalgamation basis of 
$T_{ab}^*$ if and only if $o(g)$ is a power of a prime number 
$p$, and $G$ has a unique maximal p-subgroup isomorph to 
$\mathbb{Z}/p^n\mathbb{Z}$ where $n\in\mathbb{N^*}$.
\end{lem}
\preuve 
The proof follows directly from Corollary \ref{corchar amal},
Example \ref{exm illustr} and Theorem \ref{thm char aml by uniquess}.

{\bf{Acknowledgement}}:\\
The researcher would like to thank the Deanship of Scientific Research, Qassim University.

\begin{flushleft}
Mohammed Belkasmi\\
Department of Mathematics\\  College of Sciences, Qassim University\\ P.O. Box 6644,
Buraydah 51452\\ Saudi Arabia.\\
m.belkasmi@qu.edu.sa

\end{flushleft}

 \bibliographystyle{plain}

\begin{thebibliography}{10}
\bibitem{ana}
 Belkasmi, M.
\newblock{\em Positive model theory and amalgamations}.
\newblock{\em Notre Dame Journal of Formal Logic}.Vol 55.2. (2014)
\bibitem{ana2}
Belkasmi, M.
\newblock{\em Positive amalgamations}.
\newblock{\em Logica Universalis}. Vol 14. (2020) 
\bibitem{2}
 Gobel, R and  Shelah, S.
\newblock {\em Philip Hall's problem on non-Abelian splitters}.
\newblock {\em Mathematical Proceedings of the Cambridge Philosophical Society}. Vol 134.1 (2003)
\bibitem{3}
  Karrass. M and  Solitar. D
\newblock {\em On the Free Product of Two Groups with an Amalgamated Subgroup of Finite Index in each Factor}.
\newblock {\em Proceedings of the American Mathematical Society}. Vol 26.1 (1970)
\bibitem{marker}
 Marker, D.
\newblock {\em Model Theory,  An Introduction}.
\newblock {\em Springer New York}. (2002)
\bibitem{5}
 Ordman, ET.
\newblock {\em On subgroups of amalgamated free products}.
\newblock {\em Mathematical Proceedings of the Cambridge Philosophical Society}. Vol 69.1 (1971)
\bibitem{lattice}
  Scowcroft, P.
 \newblock{\em Algebraically closed and existentially closed Abelian lattice-ordered groups}.
\newblock{\em Algebra Universalis}. (2016).
\bibitem{stable}
  Shelah, S.
 \newblock{\em Classification Theory and the Number of Non-isomorphic Models}.
\newblock{\em North-Holland Publishing Company}. (1990). 
\end{thebibliography}

\end{document}